\documentclass[letterpaper, 10pt, conference]{ieeeconf}
\IEEEoverridecommandlockouts 
\usepackage{amsmath,amssymb,url}
\usepackage{graphicx,subfigure,warmread}
\usepackage[all,import]{xy}

\newcommand{\norm}[1]{\ensuremath{\left\| #1 \right\|}}

\newcommand{\bracket}[1]{\ensuremath{\left[ #1 \right]}}

\newcommand{\parenth}[1]{\ensuremath{\left( #1 \right)}}

\newcommand{\refeqn}[1]{(\ref{eqn:#1})}
\newcommand{\reffig}[1]{Fig. \ref{fig:#1}}
\newcommand{\tr}[1]{\mbox{tr}\ensuremath{\negthickspace\bracket{#1}}}

\newcommand{\G}{\ensuremath{\mathsf{G}}}
\newcommand{\SO}{\ensuremath{\mathsf{SO(3)}}}
\newcommand{\T}{\ensuremath{\mathsf{T}}}
\renewcommand{\L}{\ensuremath{\mathsf{L}}}
\newcommand{\so}{\ensuremath{\mathfrak{so}(3)}}
\newcommand{\SE}{\ensuremath{\mathsf{SE(3)}}}

\renewcommand{\Re}{\ensuremath{\mathbb{R}}}

\newcommand{\D}{\ensuremath{\mathbf{D}}}

\newcommand{\Ad}{\ensuremath{\mathrm{Ad}}}
\newcommand{\ad}{\ensuremath{\mathrm{ad}}}
\newcommand{\g}{\ensuremath{\mathfrak{g}}}

\title{\LARGE \bf
Dynamics of Connected Rigid Bodies in a Perfect Fluid}

\author{Taeyoung Lee, Melvin Leok\authorrefmark{1}, and N. Harris McClamroch\authorrefmark{2}%
\thanks{Taeyoung Lee, Mechanical and Aerospace Engineering, Florida Institute of Technology, Melbourne, FL 39201 {\tt taeyoung@fit.edu}}%
\thanks{Melvin Leok, Mathematics, Purdue University, West Lafayette, IN 47907 {\tt mleok@math.purdue.edu}}%
\thanks{N. Harris McClamroch, Aerospace Engineering, University of Michigan, Ann Arbor, MI 48109 {\tt
nhm@umich.edu}}%
\thanks{\textsuperscript{\footnotesize\ensuremath{*}}This research has been supported in part by NSF under grants DMS-0504747 and DMS-0726263.}
\thanks{\textsuperscript{\footnotesize\ensuremath{\dagger}}This research has been supported in part by NSF under grant CMS-0555797.}
}

\begin{document}
\allowdisplaybreaks
\maketitle \thispagestyle{empty} \pagestyle{empty}

\begin{abstract}
This paper presents an analytical model and a geometric numerical integrator for a system of rigid bodies connected by ball joints, immersed in an irrotational and incompressible fluid. The rigid bodies can translate and rotate in three-dimensional space, and each joint has three rotational degrees of freedom. This model characterizes the qualitative behavior of three-dimensional fish locomotion. A geometric numerical integrator, refereed to as a Lie group variational integrator, preserves Hamiltonian structures of the presented model and its Lie group configuration manifold. These properties are illustrated by a numerical simulation for a system of three connected rigid bodies.
\end{abstract}

\section{Introduction}

Fish locomotion has been investigated in the fields of biomechanics and engineering (see \cite{Sfa.IJoOE1999} and references therein). This is a challenging problem as it involves interaction of a deformable fish body with an unsteady fluid, through which an internal muscular force of the fish is translated into an external propulsive force exerted on the fluid.

Various mathematical models of fish locomotion have been formulated. A quasi-static model based on a steady state flow theory is developed in~\cite{Tay.PotRSoLSA1952}, and an elastic plate model that treats a fish as an elongated slender body is studied in~\cite{Wu.JoFM1961,Lig.BK1975,Lig.PotRSoLSB1971}. The effects of body thickness for the slender body model are considered in~\cite{Lig.ARoFM1969}. Numerical models involving computational fluid dynamics techniques appear in~\cite{Nak.Pot4IOaPEC1994,Ram.AiFM1996}. The body of a fish is modeled as a planar articulated rigid body in~\cite{Kan.JoNS2005,Kel.1998,Rad.2003}.

The planar articulated rigid body model has become popular in engineering area, as it depicts underwater robotic vehicles that move and steer by changing their shape~\cite{Jal.PoAViMCS1995,Bar.1996}. Furthermore, if it is assumed that the ambient fluid is incompressible and irrotational, then equations of motion of the articulated rigid body can be derived without explicitly incorporating  fluid variables~\cite{Kan.JoNS2005}. The effect of the fluid is accounted by added inertia terms of the rigid body. This model is known to characterize the qualitative behavior of fish swimming properly~\cite{Kan.JoNS2005}. Based on this assumption, optimal shape changes of a planar articulated body to achieve a desired locomotion has been studied in~\cite{Kan.PotICoDaC2005,Ros.PotACC2006}.

By following~\cite{Kan.JoNS2005,Kel.1998,Rad.2003}, we consider a system of connected rigid bodies immersed in a incompressible and irrotational fluid, and we first develop an analytical model of it. The contribution of this paper is that the connected rigid bodies can freely translate and rotate in three-dimensional space, and each joint has three rotational degrees of freedom. This is important for understanding the locomotion of a fish with a blunt body and a large caudal fin.

The second part of this paper deals with a geometric numerical integrator of connected rigid bodies in a perfect fluid. Geometric numerical integration is concerned with developing numerical integrators that preserve geometric features of a system, such as invariants, symmetry, and reversibility~\cite{Hai.BK2006}. It is critical for a numerical simulation of Hamiltonian systems on a Lie group to preserve both the symplectic property of Hamiltonian flows and the Lie group structure~\cite{Lee.CMaDA2007}. A geometric numerical integrator, referred to as a Lie group variational integrator, has been developed for a Hamiltonian system on an arbitrary Lie group in~\cite{Lee.2008}.

A system of connected rigid bodies is a Hamiltonian system, and its configuration manifold is expressed as a product of the special Euclidean group and copies of the special orthogonal group. This paper develops a Lie group variational integrator for the connected rigid bodies in a perfect fluid based on the results presented in~\cite{Lee.2008}. The proposed geometric numerical integrator preserves symplecticity and momentum maps, and exhibits desirable energy properties. It also respects the Lie group structure of the configuration manifold, and avoids the singularities and complexities associated with local coordinates.

In summary, this paper develops an analytical model and a geometric numerical integrator for a system of connected rigid bodies in a perfect fluid. These provide a three-dimensional mathematical model and a reliable numerical simulation tool that characterizes the qualitative properties of fish locomotion.

This paper is organized as follows. A system of connected rigid bodies immersed in a perfect fluid is described in Section~\ref{sec:CRB}. An analytical model and a Lie group variational integrator are developed in Section~\ref{sec:AM} and in Section~\ref{sec:LGVI}, respectively, followed by a numerical example in Section~\ref{sec:NE}.

\section{Connected Rigid Bodies Immersed in a Perfect Fluid}\label{sec:CRB}

\begin{figure}
$$\begin{xy}
\xyWARMprocessEPS{RBPF3d}{eps}
\xyMarkedImport{}
\xyMarkedMathPoints{1-11}
\end{xy}$$
\caption{Connected Rigid Bodies Immersed in a Perfect Fluid}\label{fig:CRB}
\end{figure}
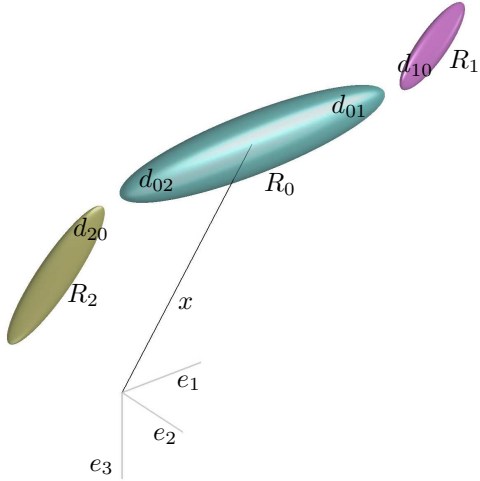

Consider three connected rigid bodies immersed in a perfect fluid. We assume that these rigid bodies are connected by a ball joint that has three rotational degrees of freedom, and the fluid is incompressible and irrotational. We also assume each body has neutral buoyancy: the mass of the body equals the mass of the fluid it displaces. This model is illustrated by \reffig{CRB}.

\renewcommand{\xyWARMinclude}[1]{\includegraphics[width=0.90\columnwidth]{#1}}

We choose a reference frame and three body-fixed frames. The origin of each body-fixed frame is located at the mass center of the rigid body and it is aligned along the principal axes. Define
{\allowdisplaybreaks
\begin{center}
\begin{tabular}{lp{6cm}}
$R_i\in\SO$ & Rotation matrix from the $i$-th body-fixed frame to the reference frame\\
$\Omega_i\in\Re^3$ & Angular velocity of the $i$-th body, represented in the $i$-th body-fixed frame\\
$x\in\Re^3$ & Vector from the origin of the reference frame to the mass center of the $0$-th body, represented in the reference frame\\
$d_{ij}\in\Re$ & Vector from the mass center of the $i$-th body to the ball joint connecting the $i$-th body with the $j$-th body, represented in the $i$-th body-fixed frame\\
$m^b_i\in\Re$ & Mass of the $i$-th body\\
$J^b_i\in\Re^{3\times3}$ & Inertia matrix of the $i$-th body
\end{tabular}
\end{center}}
\noindent for $i,j\in\{0,1,2\}$.

A configuration of this system can be described by the location of the mass center of the central body, and the attitude of each rigid body with respect to the reference frame. So, the configuration manifold is $\G=\SE\times\SO\times\SO$, where $\SO=\{R\in\Re^{3\times 3}\,|\, R^TR=I, \det{R}=1\}$, and $\SE=\SO\textcircled{s}\Re^3$.

The attitude kinematics equation is given by
\begin{align*}
    \dot R_i = R_i\hat\Omega_i
\end{align*}
for $i\in\{0,1,2\}$, where the \textit{hat map} $\hat\cdot:\Re^3\rightarrow\so$ is defined such that $\hat x y=x\times y$ for any $x,y\in\Re^3$.

\section{Continuous-time Analytical Model}\label{sec:AM}

In this section, we develop continuous-time equations of motion for a system of connected rigid bodies in a perfect fluid. As the fluid is irrotational, equations of motion can be expressed without explicitly incorporating fluid variables, and the effects of the ambient fluid is encountered by added inertia terms~\cite{Kan.JoNS2005}. To simplify expressions for the added inertia terms, we assume each body is an ellipsoid.

We first find an expression for the Lagrangian of the system, and substitute it into Euler-Lagrange equations.

\subsection{Lagrangian}
The total kinetic energy of connected rigid bodies immersed in a fluid can be written as the sum of the kinetic energy of the rigid bodies $T_{\mathcal{B}_i}$ and the kinetic energy of the fluid $T_{\mathcal{F}}$:
\begin{align*}
    T= \sum_{i=0}^2 T_{\mathcal{B}_i} + T_{\mathcal{F}}.
\end{align*}

\paragraph*{Kinetic energy of rigid bodies}

Let $V_i\in\Re^3$ be the velocity of the mass center of the $i$-th body represented in the $i$-th body-fixed frame for $i\in\{0,1,2\}$. Since $\dot x$ represents the velocity of the $0$-th rigid body in the reference frame, we obtain
\begin{align}
V_0=R_0^T\dot x.\label{eqn:V0}
\end{align}
The location of the mass center of the first rigid body can be written as $x+ R_0d_{01}-R_1d_{10}$ with respect to the reference frame. Therefore, $V_1$ is given by
\begin{align}
    V_1 & = R_1^T (\dot x + R_0\hat\Omega_0 d_{01} -R_1\hat\Omega_1 d_{10})\nonumber\\
    & = R_1^T \dot x - R_1^T R_0 \hat d_{01}\Omega_0 + \hat d_{10}\Omega_1.
\end{align}
Similarly,
\begin{align}
    V_2  & = R_2^T \dot x - R_2^T R_0 \hat d_{02}\Omega_0 + \hat d_{20}\Omega_2.\label{eqn:V2}
\end{align}

The kinetic energy of rigid bodies is given by
\begin{align}
    T_{\mathcal{B}} & = \sum_{i=0}^2 \frac{1}{2}m^b_i V_i\cdot V_i + \frac{1}{2} \Omega_i \cdot J^b_i\Omega_i.\label{eqn:TB}
\end{align}

\begin{figure*}[!t]\footnotesize\selectfont
\begin{align}
\mathbb{I}=
    \begin{bmatrix}
    J_0-\hat d_{01}R_0^TR_1M_1R_1^T R_0\hat d_{01}-\hat d_{02}R_0^T R_2M_2R_2^T R_0\hat d_{02}
    & \hat d_{01}R_0^TR_1M_1R_1^T+\hat d_{02}R_0^TR_2M_2R_2^T
    & \hat d_{01} R_0^T R_1M_1 \hat d_{10}
    & \hat d_{02}R_0^T R_2 M_2\hat d_{20}\\
    -R_1M_1R_1^TR_0\hat d_{01}-R_2M_2R_2^T R_0\hat d_{02}
    & R_0M_0R_0^T+R_1M_1R_1^T+R_2M_2R_2^T
    & R_1M_1\hat d_{10}
    & R_2M_2\hat d_{20}\\
    \hat d_{10}M_1R_1^T R_0 \hat d_{01}
    & -\hat d_{10}M_1R_1^T
    & J_1-\hat d_{10}M_1\hat d_{10}
    & 0\\
    \hat d_{20}M_2R_2^T R_0\hat d_{02}
    & -\hat d_{20} M_2 R_2^T
    & 0
    & J_2-\hat d_{20}M_2\hat d_{20}
    \end{bmatrix}\label{eqn:II}
\end{align}
\hrulefill
\end{figure*}

\paragraph*{Kinetic energy of fluid}

The kinetic energy of the fluid is given by
\begin{align*}
    T_{\mathcal{F}}= \frac{1}{2}\int_{\mathcal{F}} \rho_f \norm{u}^2 dv,
\end{align*}
where $\rho_f$ is the density of the fluid, $u$ is the velocity field of the fluid and $dv$ is the standard volume element in $\Re^3$. Since the flow is irrotational, the velocity field can be expressed as a gradient of a potential. Under these conditions, the kinetic energy of the fluid can be written as
\begin{align*}
    T_{\mathcal{F}} & = \sum_{i,j=0}^2 \frac{1}{2}M^f_{ij} V_i\cdot V_j + \frac{1}{2} \Omega_i \cdot J^f_{ij} \Omega_j+ D^f_{ij} V_i\cdot \Omega_j,
\end{align*}
where $M^f_{ij}, J^f_{ij}, D^f_{ij} \in\Re^{3\times 3}$ are referred to as \textit{added inertia matrices}~\cite{Hol.PD1998}. Here we assume that the flow near one rigid body is not affected by other rigid bodies: the added inertia matrices $M^f_{ij}, J^f_{ij}, D^f_{ij}$ are equal to zero when $i\neq j$. The resulting model captures the qualitative properties of the interaction between rigid body dynamics and fluid dynamics correctly \cite{Kan.JoNS2005,Kan.PotICoDaC2005}.

Expressions for added inertia matrices for an ellipsoidal body are derived in \cite{Lam.BK1932}. Let $l_q\in\Re$ be the length of the $q$-th principal axis of an ellipsoid for $q\in \{1,2,3\}$. Define constants
\begin{align*}
    \gamma_q = l_1l_2l_3 \int_{0}^\infty \frac{d\nu}{(l_q^2+\nu)\sqrt{(l_1^2+\nu)(l_2^2+\nu)(l_3^2+\nu)}}
\end{align*}
for $q\in \{1,2,3\}$ and
\begin{align*}
    \lambda_1 = \frac{1}{5}m^b \parenth{\frac{(l_2^2-l_3^2)^2(\gamma_3-\gamma_2)}{2(l_2^2-l_3^2)+(l_2^2+l_3^2)(\gamma_2-\gamma_3)}}.
\end{align*}
Constants $\lambda_2$ and $\lambda_3$ are given by cyclic permutations of this expression. Then, the added inertia matrices of the ellipsoid are given by
\begin{align}
    M^f & = m^b\mathrm{diag}\bracket{\frac{\gamma_1}{2-\gamma_1},\,\frac{\gamma_2}{2-\gamma_2},\,\frac{\gamma_3}{2-\gamma_3}},\\
    J^f & = \mathrm{diag}\bracket{\lambda_1,\,\lambda_2,\,\lambda_3},\\
    D^f & = 0.
\end{align}
Using these expressions, we find added inertia matrices $M^f_{ii},J^f_{ii}$ for each rigid body.

In summary, the kinetic energy of the fluid surrounding ellipsoidal rigid bodies is given by
\begin{align}
    T_{\mathcal{F}}= \sum_{i=0}^2 \frac{1}{2} M^f_{ii} V_i\cdot V_i + \frac{1}{2} \Omega_i \cdot J^f_{ii} \Omega_i.\label{eqn:TF}
\end{align}

\paragraph*{Total kinetic energy}

Define total inertia matrices
\begin{align}
    M_i & = m^b_iI_{3\times 3}  + M^f_{ii},\\
    J_i & = J^b_i + J^f_{ii}
\end{align}
for $i=\{0,1,2\}$. From \refeqn{TB} and \refeqn{TF}, the total kinetic energy is given by
\begin{align}
    T & = \sum_{i=0}^2 \frac{1}{2}M_i V_i\cdot V_i + \frac{1}{2} \Omega_i \cdot J_i\Omega_i.
\end{align}
Substituting \refeqn{V0}-\refeqn{V2}, this can be written as
\begin{align}
T = \frac{1}{2} \xi^T \mathbb{I}(R_0,R_1,R_2) \xi,\label{eqn:T}
\end{align}
where $\xi=[\Omega_0;\dot x;\Omega_1;\Omega_2]\in\Re^{12}$ and the matrix $\mathbb{I}(R_0,R_1,R_2)\in\Re^{12\times 12}$ is given by \refeqn{II}. Since there is no potential field, this is equal to the Lagrangian of the connected rigid bodies immersed in a perfect fluid.

\subsection{Euler-Lagrange Equations}

Euler-Lagrange equations for a mechanical system that evolves on an arbitrary Lie group are given by
\begin{gather}
\frac{d}{dt}\D_\xi L(g,\xi)-\ad^*_\xi \cdot \D_\xi L(g,\xi) -\T_e^*\L_g\cdot \D_g L(g,\xi)=0,\label{eqn:EL0}\\
\dot g= g\xi,\label{eqn:EL1}
\end{gather}
where $L:\T\G\simeq \G\times\g\rightarrow \Re$ is the Lagrangian of the system~\cite{Lee.2008}. Here $\D_\xi L(g,\xi)\in\g^*$ denotes the derivative of the Lagrangian with respect to $\xi\in\g$, $\ad^*:\g\times\g^*\rightarrow\g^*$ is $\mathrm{co}$-$\mathrm{ad}$ operator, and $\T_e^*\L_g:\T^*\G\rightarrow\g^*$ denotes the cotangent lift of the left translation map $\L_g:\G\rightarrow\G$ (see \cite{Mar.BK1999} for the detailed definitions).

Using this result, we develop Euler-Lagrange equations of a system of connected rigid bodies in a perfect fluid. To simplify the derivation, we consider the configuration manifold given by $\G=\SO\times\Re^3\times\SO\times\SO$, left-trivialize $\T\G$ to yield $\G\times\g$, and identify its Lie algebra $\g$ with  $\Re^{12}$ by the hat map. For $\xi=[\Omega_0;\dot x;\Omega_1;\Omega_2]\in\g$ and $p=[p_0;p_x;p_1;p_2]\in\g^*$, the $\mathrm{co}$-$\mathrm{ad}$ operator is given by $\ad^*_\xi p =[-\hat\Omega_0p_0;p_x;-\hat\Omega_1 p_1;-\hat\Omega_2 p_2]$.

\paragraph*{Derivatives of the Lagrangian} The derivative of the Lagrangian with respect $\xi$ is given by
\begin{align}
    \D_\xi L (g,\xi) = \mathbb{I}(R_0,R_1,R_2)\xi.\label{eqn:DxiL}
\end{align}
The derivative of the Lagrangian with respect to $g=(R_0,x,R_1,R_2)\in\G$ can be written as
\begin{align}
    &\T_e^*\L_g\cdot \D_g L(g,\xi)\nonumber\\
    &=[\T_I^*\L_{R_0}\cdot \D_{R_0} L;\, \D_x L;\,
    \T_I^*\L_{R_1}\cdot \D_{R_1} L;\, \T_I^*\L_{R_2}\cdot \D_{R_2} L].
\end{align}
An expression for the first term of this can be found as follows. For any $\eta_0\in\Re^3$, let $g_0^\epsilon=[R_0\exp\epsilon\eta_0,x,R_1,R_2]\in\G$. Then, we have
\begin{align*}
&(\T_I^* \L_{R_0} \cdot \D_{R_0}L)\cdot \eta_0  =
\frac{d}{d\epsilon}\bigg|_{\epsilon=0} L(g_0^\epsilon,\xi)
\\
&=-\dot x^T R_0 M_0\hat \eta_0 R_0^T\dot x
+ \sum_{i=1}^2 \Big( -\Omega_0^T \hat d_{0i}R_0^T R_i M_i R_i^T R_0 \hat \eta_0\hat d_{0i} \Omega_0\\
& \quad - \dot x^T R_i M_i R_i^T R_0 \hat\eta_0 \hat d_{0i} \Omega_0
- \Omega_0^T \hat d_{0i} \hat\eta_0 R_0^T R_i M_i \hat d_{i0} \Omega_i\Big)\nonumber\\
& = \parenth{-\widehat{R_0^T\dot x}M_0R_0^T\dot x
-\sum_{i=1}^2
 \widehat{\hat d_{0i} \Omega_0} R_0^T R_i M_i V_i }\cdot \eta_0,
\end{align*}
where we use identities: $x\cdot y =x^Ty=y^Tx$, $\hat x y =-\hat y x$ for any $x,y\in\Re^3$. Since this is satisfied for any $\eta_0\in\Re^3$, we obtain
\begin{align}
    \T_I^* \L_{R_0} \cdot \D_{R_0}L = -\widehat{R_0^T\dot x}M_0R_0^T\dot x
-\sum_{i=1}^2
 \widehat{\hat d_{0i} \Omega_0} R_0^T R_i M_i V_i.
\end{align}
Similarly, we find
\begin{align}
\D_x L & = 0,\\
\T_I^* \L_{R_i} \cdot \D_{R_i}L & = \widehat{M_iV_i} R_i^T (\dot x - R_0\hat d_{0i}\Omega_0)\label{eqn:DRiL}
\end{align}
for $i\in\{1,2\}$.

\paragraph*{Euler-Lagrange Equations}
Substituting \refeqn{DxiL}-\refeqn{DRiL} into \refeqn{EL0}-\refeqn{EL1}, and rearranging, Euler-Lagrange equations for the connected rigid bodies immersed in a perfect fluid are given by
\begin{gather}
\begin{aligned}
&\mathbb{I}(R_0,R_1,R_2)
    \begin{bmatrix} \dot\Omega_0 \\ \ddot x \\ \dot\Omega_1 \\ \dot\Omega_2\end{bmatrix}\\
&+\begin{bmatrix}
    \Omega_0\times J\Omega_0 +\widehat{ R_0^T\dot x} M_0 R_0^T \dot x
    + \sum_{i=1}^2\hat d_{0i}R_0^T R_i W_i \\
    R_0(\hat\Omega_0M_0-M_0\hat\Omega_0)R_0^T\dot x +\sum_{i=1}^2 R_i W_i\\
    \Omega_1\times J_1\Omega_1 +V_1\times M_1V_1 -\hat d_{10}W_1\\
    \Omega_2\times J_2\Omega_2 +V_2\times M_2 V_2 -\hat d_{20}W_2
\end{bmatrix}=0,\end{aligned}\\
\dot R_0=R\hat\Omega_1,\quad \dot R_1=R_1\hat\Omega_1,\quad \dot R_2=R_2\hat\Omega_2,
\end{gather}
where
\begin{align}
V_i & = R_i^T \dot x - R_i^T R_0 \hat d_{0i} \Omega_0 + \hat d_{i0}\Omega_i,\\
W_i & =(\hat\Omega_iM_i-M_i\hat\Omega_i)(R_i^T\dot x - R_i^TR_0\hat d_{0i}\Omega_0)\nonumber\\
&\quad -M_iR_i^TR_0 \hat\Omega_0 \hat d_{0i}\Omega_0 + \hat\Omega_i M_i \hat d_{i0}\Omega_i
\end{align}
for $i\in\{1,2\}$.

\paragraph*{Hamilton's equations}

Let the momentum of the system be $\mu=[p_0;p_x;p_1;p_2]\in\Re^{12}\simeq \g^*$. The Legendre transformation is given by
\begin{align}
    \mu = \D_\xi L(g,\xi) = \mathbb{I}(R_0,R_1,R_2)\xi.
\end{align}
The corresponding Hamilton's equations can be written as
\begin{align}
\dot p_0 &= -\hat\Omega_0 p_0-\widehat{R_0^T\dot x}M_0R_0^T\dot x
-\sum_{i=1}^2
 \widehat{\hat d_{0i} \Omega_0} R_0^T R_i M_i V_i,\label{eqn:dotp0}\\
\dot p_x &= 0,\\
\dot p_i &= -\hat\Omega_i p_i + \widehat{M_iV_i} R_i^T (\dot x - R_0\hat d_{0i}\Omega_0)\label{eqn:dotpi}
\end{align}
for $i\in\{1,2\}$.

\paragraph*{Conserved quantities} As the Lagrangian is invariant under rigid translation and rotation of the entire system, the total linear momentum $p_x\in\Re^3$ and the total angular momentum $p_\Omega=\hat x p_x +\sum_{i=0}^2 R_i p_i \in\Re^3$ are preserved.

\section{Lie Group Variational Integrator}\label{sec:LGVI}

The continuous-time Euler-Lagrange equations and Hamilton's equations developed in the previous section provide analytical models of the connected rigid bodies in a perfect fluid. However, they are not suitable for a numerical study since a direct numerical integration of those equations using a general purpose numerical integrator, such as an explicit Runge Kutta method, may not preserve the geometric properties of the system accurately~\cite{Hai.BK2006}.

Variational integrators provide a systematic method of developing geometric numerical integrators for Lagrangian/Hamiltonian systems~\cite{Mar.2001}. As it is derived from a discrete analogue of Hamilton's principle, it preserves symplecticity and the momentum map, and it exhibits good total energy behavior. Lie group methods conserve the structure of a Lie group configuration manifold as it updates a group element using the group operation~\cite{Ise.2000}.

These two methods have been unified to obtain a Lie group variational integrator for Lagrangian/Hamiltonian systems evolving on a Lie group~\cite{Lee.2008}. This preserves symplecticity and group structure of those systems concurrently. It has been shown that this property is critical for accurate and efficient simulations of rigid body dynamics~\cite{Lee.CMaDA2007}.

In this section, we develop a Lie group variational integrator for the connected rigid bodies in a perfect fluid. We first obtain an expression for a discrete Lagrangian and substitute it into the discrete-time Euler-Lagrange equations.

\subsection{Discrete Lagrangian}

Let $h>0$ be a fixed integration step size, and let a subscript $k$ denote the value of a variable at the $k$-th time step. We define a discrete-time kinematics equation as follows. Define $f_k=(F_{0_k},\Delta x_k,F_{1_k},F_{2_k})\in\G$ for $\Delta x_k\in\Re^3$, $F_{0_k},F_{1_k},F_{2_k}\in\SO$ such that $g_{k+1}=g_k f_k$:
\begin{align}
    (R_{0_{k+1}},\,& x_{k+1},\,R_{1_{k+1}},\,R_{2_{k+1}})\nonumber\\
    & =(R_{0_{k}}F_{0_k},\,x_{k}+\Delta x_k,\,R_{1_{k}}F_{1_k},\,R_{2_{k}}F_{2_k}).
\end{align}
Therefore, $f_k$ represents the relative update between two integration steps. This ensures that the structure of the Lie group configuration manifold is numerically preserved.

A discrete Lagrangian $L_d(g_k,f_k):\G\times\G\rightarrow\Re$ is an approximation of the Jacobi solution of the Hamilton--Jacobi equation, which is given by the integral of the Lagrangian along the exact solution of the Euler-Lagrange equations over a single time step:
\begin{align*}
    L_d(g_k,f_k)\approx \int_0^h L(\tilde g(t),{\tilde g}^{-1}(t)\dot{\tilde g} (t))\,dt,
\end{align*}
where $\tilde g(t):[0,h]\rightarrow \G$ satisfies Euler-Lagrange equations with boundary conditions $\tilde{g}(0)=g_k$, $\tilde{g}(h)=g_kf_k$. The resulting discrete-time Lagrangian system, referred to as a variational integrator, approximates the Euler-Lagrange equations to the same order of accuracy as the discrete Lagrangian approximates the Jacobi solution.

The kinetic energy given by \refeqn{T} can be rewritten as
\begin{align*}
T   & = \frac{1}{2} \dot x^TR_0 M_0 R_0^T\dot x+ \frac{1}{2} \Omega_0^T J_0\Omega_0
\\&\quad+\sum_{i=1}^2 \Big( \frac{1}{2}\dot x^T R_iM_iR_i^T
+ \frac{1}{2} \Omega_i^T (J_i-\hat d_{i0}M_i\hat d_{i0})\Omega_i\nonumber\\
&\quad\quad  +\frac{1}{2}d_{0i}^T\dot R_0^T R_i M_i R_i^T \dot R_0 d_{0i}
 + \dot x^T R_i M_i R_i^T \dot R_0 d_{0i}
\\&\quad\quad - \dot x^T R_iM_i\hat\Omega_i d_{i0}
-  d_{0i}^T \dot R_0^T R_iM_i \hat\Omega_i d_{i0}\Big).
\end{align*}
>From this, we choose the discrete Lagrangian as
\begin{align}
    L_{d_k} & = \frac{1}{2h}\Delta x_k^T R_{0_k}M_0 R_{0_k}^T\Delta x_k +\frac{1}{h}\tr{(I-F_{0_k})  J_{d_0}}\nonumber\\
    &\quad  +\sum_{i=1}^2 \Big( \frac{1}{2h}\Delta x_k^T
    R_{i_k}M_i R_{i_k}^T\Delta x_k + \frac{1}{h}\tr{ (I-F_{i_k}) J_{d_i}'}\nonumber\\
    &\quad + \frac{1}{2h} d_{0i}^T(F_{0_k}^T-I)R_{0_k}^T R_{i_k} M_i R_{i_k}^T R_{0_k} (F_{0_k}-I)d_{0i}\nonumber\\
    &\quad +\frac{1}{h}\Delta x_k^T R_{i_k} M_i R_{i_k}^T R_{0_k}(F_{0_k}-I)d_{0i}\nonumber\\
    &\quad -\frac{1}{h}\Delta x_k^T R_{i_k} M_i (F_{i_k}-I)d_{i0}\nonumber\\
    &\quad-\frac{1}{h}d_{0i}^T (F_{0_k}^T -I) R_{0_k}^T R_{i_k} M_i (F_{i_k}-I) d_{i0}\Big),\label{Ldk}
\end{align}
where nonstandard inertia matrices are defined as
\begin{gather}
    J_{d_0} = \frac{1}{2}\tr{J_0}I-J_0,\label{eqn:Jd0}\\
    J_{d_i}' = \frac{1}{2}\tr{J_i'}I-J_i',\quad
    J_i' = J_i - \hat d_{i0}M_i \hat d_{i0},\label{eqn:Jdi}
\end{gather}
for $i\in\{1,2\}$.

\subsection{Discrete-time Euler-Lagrange Equations}

For a discrete Lagrangian on $\G\times\G$, the following discrete-time Euler-Lagrange equations, referred to as a Lie group variational integrator, were developed in~\cite{Lee.2008}.
\begin{gather}
\begin{aligned}
    \T_e\L_{f_{k}}\cdot \D_{f_k}L_{d_k}-&\Ad^*_{f_{k+1}^{-1}}\cdot(\T_e\L_{f_{k+1}}\cdot \D_{f_{k+1}}L_{d_{k+1}})\\
    &\quad+\T_e\L_{g_{k+1}}\cdot \D_{g_{k+1}} L_{d_{k+1}}=0,
\end{aligned}\label{eqn:DEL0}\\
    g_{k+1} = g_k f_k,\label{eqn:DEL1}
\end{gather}
where $\Ad^*:\G\times\g^*\rightarrow\g^*$ is $\mathrm{co}$-$\mathrm{Ad}$ operator~\cite{Mar.BK1999}.

Using this result, we develop a Lie group variational integrator for connected rigid bodies in a perfect fluid. For $f=(F_0,\Delta x,F_1,F_2)\in \G$ and $p=[p_0;p_x;p_1;p_2]\in\g^*\simeq\Re^{12}$, the $\mathrm{co}$-$\mathrm{Ad}$ operator is given by $\Ad^*_{f^{-1}} p = [F_0p_0;p_x;F_1p_1;F_2p_2]=[(F_0\hat p_0 F_0^T)^\vee;p_x;(F_1\hat p_1 F_1^T)^\vee;(F_2\hat p_2 F_2^T)^\vee]$, where the \textit{vee map} $\vee:\so\rightarrow\Re^3$ denotes the inverse of the hat map.

\paragraph*{Derivatives of the discrete Lagrangian}
We find expressions for the derivatives of the discrete Lagrangian. The derivative of the discrete Lagrangian with respect to $F_{0_k}$ is given by
\begin{align*}
\D_{F_{0_k}}L_{d_k}\cdot \delta F_{0_k} & = \frac{1}{h}\tr{-\delta F_{0_k} J_{d_0}}
+\frac{1}{h}\sum_{i=1}^2 A_{i_k}^T R_{0_k}\delta F_{0_k} d_{0i},
\end{align*}
where we define, for $i\in\{1,2\}$,
\begin{align}
A_{i_k} &= R_{i_k} M_i \parenth{R_{i_k}^T  B_{i_k} - (F_{i_k}-I)d_{i0}},\label{eqn:Aik}\\
B_{i_k} & = \Delta x_k + R_{0_k} (F_{0_k}-I)d_{0i}.\label{eqn:Bik}
\end{align}
The variation of $F_{0_k}$ can be written as $\delta F_{0_k}=F_{0_k}\hat\zeta_{0_k}$ for $\zeta_{0_k}\in\Re^3$. Therefore, we have
\begin{align*}
\D_{F_{0_k}}L_{d_k}& \cdot (F_{0_k}\hat\zeta_{0_k})
 =  (\T_I^* \L_{F_{0_k}}\cdot \D_{F_{0_k}}L_{d_k}) \cdot \zeta_{0_k}\\
& = \frac{1}{h}\tr{- F_{0_k} \hat\zeta_{0_k} J_{d_0}}
+\frac{1}{h}\sum_{i=1}^2 A_{i_k}^T R_{0_{k+1}} \hat\zeta_{0_k} d_{0i}.
\end{align*}
By repeatedly applying a property of the trace operator, $\mbox{tr}[AB]=\mbox{tr}[BA]=\mbox{tr}[A^TB^T]$ for any $A,B\in\Re^{3\times 3}$, the first term can be written as $\mbox{tr}\{-F_{0_k} \hat\zeta_{0_k} J_{d_0}\}= \mbox{tr}\{- \hat\zeta_{0_k} J_{d_0}F_{0_k}\}=\mbox{tr}\{\hat\zeta_{0_k}F_{0_k}^TJ_{d_0}\}=-\frac{1}{2}\mbox{tr}\{\hat\zeta_{0_k}(J_{d_0}F_{0_k}-F_{0_k}^T J_{d_0})\}$. Using a property of the hat map, $x^T y = -\frac{1}{2}\mbox{tr}[\hat x\hat y]$ for any $x,y\in\Re^3$, this can be further written as $((J_{d_0}F_{0_k}-F_{0_k}^T J_{d_0})^\vee) \cdot \zeta_{0_k}$. As $\hat x y =-\hat y x$ for any $x,y\in\Re^3$, the second term can be written as $A_{i_k}^T R_{0_{k+1}} \hat\zeta_{0_k} d_{0i}=-A_{i_k}^T R_{0_{k+1}} \hat d_{0i} \zeta_{0_k}
=(\hat d_{0i} R_{0_{k+1}}^T A_{i_k})\cdot \zeta_{0_k}$. Using these, we obtain
\begin{align}
\T_I^* & \L_{F_{0_k}}\cdot \D_{F_{0_k}}L_{d_k} \nonumber\\
& = \frac{1}{h}(J_{d_0}F_{0_k}-F_{0_k}^T J_{d_0})^\vee+\frac{1}{h}\sum_{i=1^2}\hat d_{0i} R_{0_{k+1}}^T A_{i_k}.\label{eqn:DF0Ld}
\end{align}
Similarly, we can derive the derivatives of the discrete Lagrangian as follows.
\begin{align}
&\T^*_I\L_{F_{i_k}}\cdot \D_{F_{i_k}}L_{d_k} = \nonumber\\
&\quad \frac{1}{h} (J'_{i_d}F_{i_k}-F_{i_k}^T J'_{i_d})^\vee
  -\frac{1}{h}\hat d_{i0} F_{i_k}^T M_i R_{i_k}^TB_{i_k},\\
& \D_{\Delta x_k} L_{d_k} =\frac{1}{h} R_{0_k}M_0 R_{0_k}^T \Delta x_k + \frac{1}{h} A_{1_k} + \frac{1}{h} A_{2_k},\\
&\T^*_I \L_{R_{0_k}}\cdot \D_{R_{0_k}}L_{d_k} =\nonumber\\
&\quad \frac{1}{h} (M_0R_{0_k}^T \Delta x_k)^\wedge R_{0_k}^T\Delta x_k
 +\frac{1}{h}\sum_{i=1}^2 ((F_{0_k}-I)d_{0i})^\wedge R_{0_k}^TA_{i_k},\\
&\T^*_I\L_{R_{i_k}}\cdot \D_{R_{i_k}}L_{d_k} =\frac{1}{h} R_{i_k}^T\hat A_{i_k}B_{i_k}.
\label{eqn:DRiLd}
\end{align}

\paragraph*{Discrete-time Euler-Lagrange Equations}
Substituting \refeqn{DF0Ld}--\refeqn{DRiLd} into \refeqn{DEL0}-\refeqn{DEL1}, and rearranging, discrete-time Euler-Lagrange equations for the connected rigid bodies immersed in a perfect fluid are given by
\begin{gather}
\begin{aligned}
(J_{0_d}F_{0_k}- F_{0_k}^T&J_{0_d})^\vee-(F_{0_{k+1}}J_{0_d}- J_{0_d}F_{0_{k+1}}^T)^\vee\\
&+(M_0 R_{0_{k+1}}^T \Delta x_{k+1})^\wedge R_{0_{k+1}}^T\Delta x_{k+1}\\
&+ \sum_{i=1}^2 \hat d_{0i} R_{0_{k+1}}^T (A_{i_k}- A_{i_{k+1}})=0,
\end{aligned}\label{eqn:F0kp}\\
\begin{aligned}
(J'_{i_d}&F_{i_k}- F_{i_k}^TJ'_{i_d})^\vee
-(F_{i_{k+1}}J'_{i_d}- J'_{i_d}F_{i_{k+1}}^T)^\vee\\
&- \hat d_{i0}F_{i_k}^T  M_i R_{i_k}^TB_{i_k}\\
&+(\widehat {F_{i_{k+1}} d_{i0}}  M_i R_{i_{k+1}}^T+R_{i_{k+1}}^T \hat A_{i_{k+1}})B_{i_{k+1}}=0,
\end{aligned}\\
\begin{aligned}
R_{0_k}&M_0 R_{0_k}^T \Delta x_k + A_{1_k} + A_{2_k}\\
&-R_{0_{k+1}}M_0 R_{0_{k+1}}^T \Delta x_{k+1} - A_{1_{k+1}} - A_{2_{k+1}}=0,\\
\end{aligned}\label{eqn:delxkp}\\
R_{0_{k+1}}=R_{0_k} F_{0_k},\label{eqn:R0kp}\\
R_{i_{k+1}}=R_{i_k} F_{i_k},\\
x_{k+1} = x_k + \Delta x_k,\label{eqn:xkp}
\end{gather}
where inertia matrices are given by \refeqn{Jd0}, \refeqn{Jdi}, and $A_{i_k},B_{i_k}\in\Re^{3}$ are given by \refeqn{Aik}, \refeqn{Bik} for $i\in\{1,2\}$. For given $(g_0,f_0)\in\G\times\G$, $g_1\in\G$ is obtained by \refeqn{R0kp}--\refeqn{xkp}, and $f_1\in\G$ is obtained by solving \refeqn{F0kp}--\refeqn{delxkp}. This yields a discrete-time Lagrangian flow map $(g_0,f_0)\rightarrow(g_1,f_1)$, and this process is repeated.

\paragraph*{Discrete-time Hamilton's Equations}
Discrete-time Legendre transformation is given by
\begin{align*}
    \mu_k & = -\T_e^*\L_{g_k} \cdot \D_{g_k}L_{d_k} +\Ad^*_{f_k^{-1}}\cdot (\T_e^*\L_{f_k} \cdot \D_{f_k}L_{d_k}).
\end{align*}
Substituting this into discrete-time Euler-Lagrange equations, we obtain discrete-time Hamilton's equations for the connected rigid bodies immersed in a perfect fluid as follows.
\begin{align}
h p_{0_k} & = (F_{0_k}J_{0_d}- J_{0_d}F_{0_k}^T)^\vee
 -(M_0R_{0_k}^T \Delta x_k)^\wedge R_{0_k}^T\Delta x_k\nonumber\\
&\quad+\sum_{i=1}^2 \hat d_{0i} R_{0_k}^TA_{i_k},\label{eqn:p0k}\\
h p_{i_k} & = (F_{i_k}J'_{i_d}- J'_{i_d}F_{i_k}^T)^\vee
-\frac{1}{h} \widehat {F_{i_k} d_{i0}}  M_i R_{i_k}^TB_{i_k}\nonumber\\
&\quad-\frac{1}{h} R_{i_k}^T\hat A_{i_k}B_{i_k},\\
h p_{x_k} & = R_{0_k}M_0 R_{0_k}^T\Delta x_k+ A_{1_k} + A_{2_k},\label{eqn:pxk}\\
R_{0_{k+1}}&=R_{0_k} F_{0_k},\label{eqn:R0kpH}\\
R_{i_{k+1}}&=R_{i_k} F_{i_k},\\
x_{k+1}& =x_k+\Delta x_k,\label{eqn:xkpH}\\
h p_{0_{k+1}} & = (J_{0_d}F_{0_k}- F_{0_k}^TJ_{0_d})^\vee
+ \sum_{i=1}^2 \hat d_{0i} R_{0_{k+1}}^T A_{i_k},\label{eqn:p0kp}\\
h p_{i_{k+1}} & =  (J'_{i_d}F_{i_k}- F_{i_k}^TJ'_{i_d})^\vee
- \hat d_{i0}F_{i_k}^T  M_i R_{i_k}^TB_{i_k},\\
p_{x_{k+1}} & = p_{x_k}\label{eqn:pxkp},
\end{align}
where inertia matrices are given by \refeqn{Jd0}, \refeqn{Jdi}, and $A_{i_k},B_{i_k}\in\Re^{3}$ are given by \refeqn{Aik}, \refeqn{Bik} for $i\in\{1,2\}$. For given $(g_0,\mu_0)\in\G\times\g^*$, $f_1\in\G$ is obtained by solving \refeqn{p0k}--\refeqn{pxk}, and $g_1\in\G$ is given by \refeqn{R0kpH}--\refeqn{xkpH}. The momenta at the next step is obtained by \refeqn{p0kp}--\refeqn{pxkp}. This yields a discrete-time Hamiltonian flow map $(g_0,\mu_0)\rightarrow(g_1,\mu_1)$, and this process is repeated.

\section{Numerical Example}\label{sec:NE}

We show computational properties of the Lie group variational integrator developed in the previous section. The principal axes of each ellipsoid are given by
\begin{align*}
\text{Body 0: } & l_1=8,\quad l_2 =1.5,\quad l_3=2\;(\mathrm{m}),\\
\text{Body 1,2: } & l_1=5,\quad l_2 =0.8,\quad l_3=1.5\;(\mathrm{m}).
\end{align*}
We assume the density of fluid is $\rho=1\mathrm{kg/m^3}$. The corresponding inertia matrices are given by
\begin{gather*}
    M_0= \mathrm{diag}[1.0659,\,2.1696,\,1.6641],\;(\mathrm{kg})\\
    M_1=M_2=\mathrm{diag}[0.2664,\,    0.6551,\,    0.3677]\;(\mathrm{kg}),\\
    J_0= \mathrm{diag}[1.3480,\,   20.1500,\,   25.3276]\;(\mathrm{kgm^2}),\\
    J_1=J_2=\mathrm{diag}[0.1961,\,    1.7889,\,    2.9210]\;(\mathrm{kgm^2}).
\end{gather*}
The location of the ball joints with respect to the mass center of each body are chosen as
\begin{gather*}
d_{01}=-d_{02}=[8.8,\,0,\,0],\quad d_{10}=-d_{20}=[5.5,\,0,\,0]\;(\mathrm{m}).
\end{gather*}

The initial conditions are as follows:
\begin{gather*}
    R_{0_0}=I,\quad \Omega_{0_0}=[0.2,\,0.1,\,0.5]\;(\mathrm{rad/s}),\\
    R_{1_0}=I,\quad \Omega_{1_0}=[0.1,\,-0.3,\,-0.2]\;(\mathrm{rad/s}),\\
    R_{2_0}=I,\quad \Omega_{2_0}=[-0.1,\,0.4,\,-0.6]\;(\mathrm{rad/s}),\\
    x_0=[0,\,0,\,0]\;(\mathrm{m}),\quad \dot x_0=[0,\,   -0.4142,\,   -0.5900]\;(\mathrm{m/s}).
\end{gather*}
The corresponding total linear momentum is zero. These initial conditions provide a nontrivial rotational maneuver of the connected rigid bodies (an animation illustrating this maneuver is available at {\small\selectfont\url{http://my.fit.edu/~taeyoung}}).

We compute discrete-time Hamiltonian flow according to \refeqn{p0k}--\refeqn{pxkp}, and as comparison, we numerically integrate the continuous-time Hamilton's equations \refeqn{dotp0}--\refeqn{dotpi} using an explicit, variable step size, Runge-Kutta method. The timestep of the Lie group variational integrator is $h=0.001$ and the maneuver time is $100$ seconds.

\begin{figure*}
\centerline{
    \subfigure[Angular Velocity of Body 0, $\Omega_0$]{
        \includegraphics[width=0.95\columnwidth]{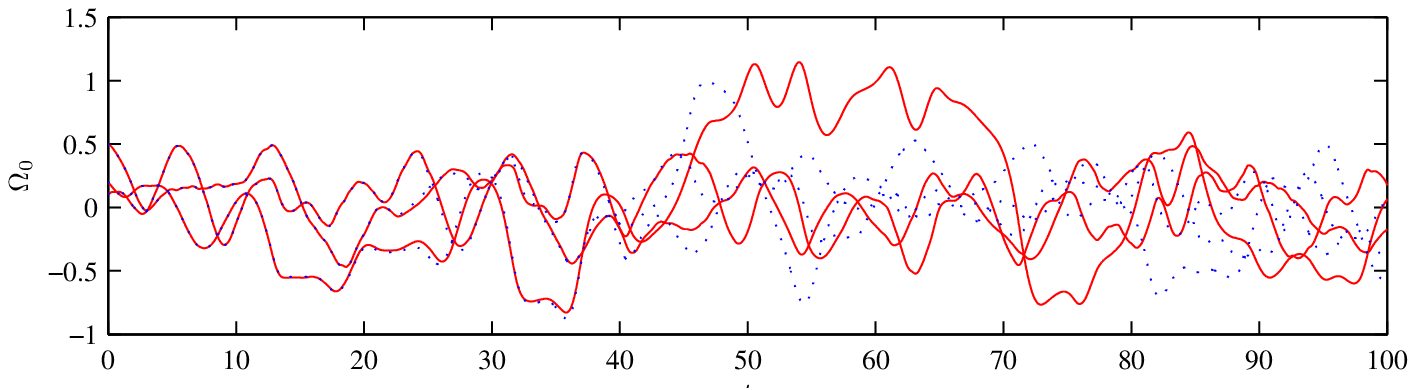}}
    \hspace*{0.05\columnwidth}
    \subfigure[Total Energy]{
        \includegraphics[width=0.95\columnwidth]{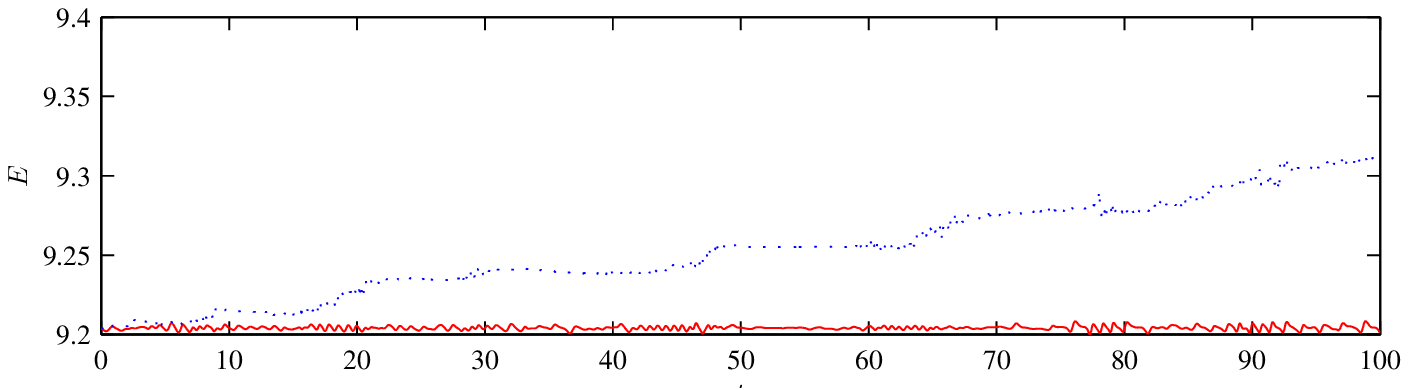}\label{fig:E}}
}
\centerline{
    \subfigure[Angular Velocity of Body 1, $\Omega_1$]{
        \includegraphics[width=0.95\columnwidth]{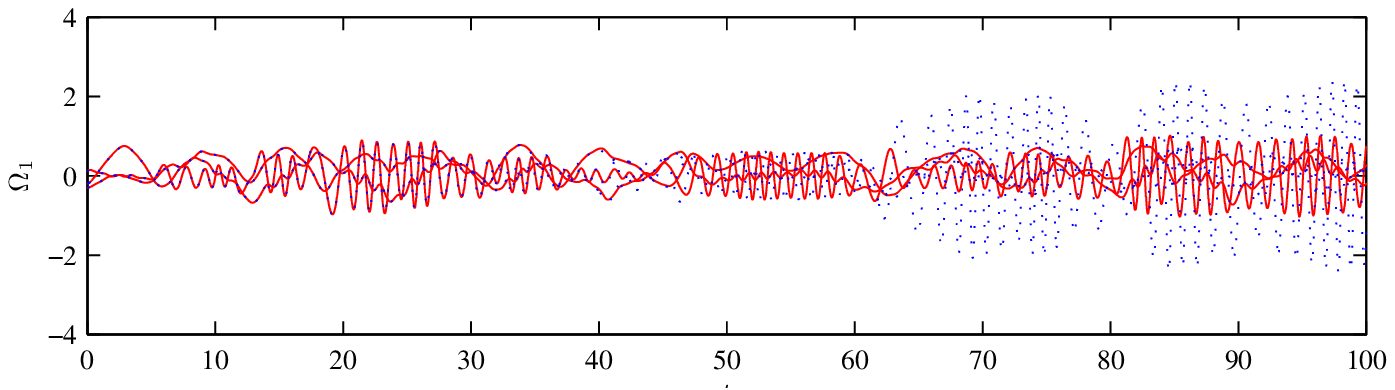}}
    \hspace*{0.05\columnwidth}
    \subfigure[Total Linear Momentum]{
        \includegraphics[width=0.95\columnwidth]{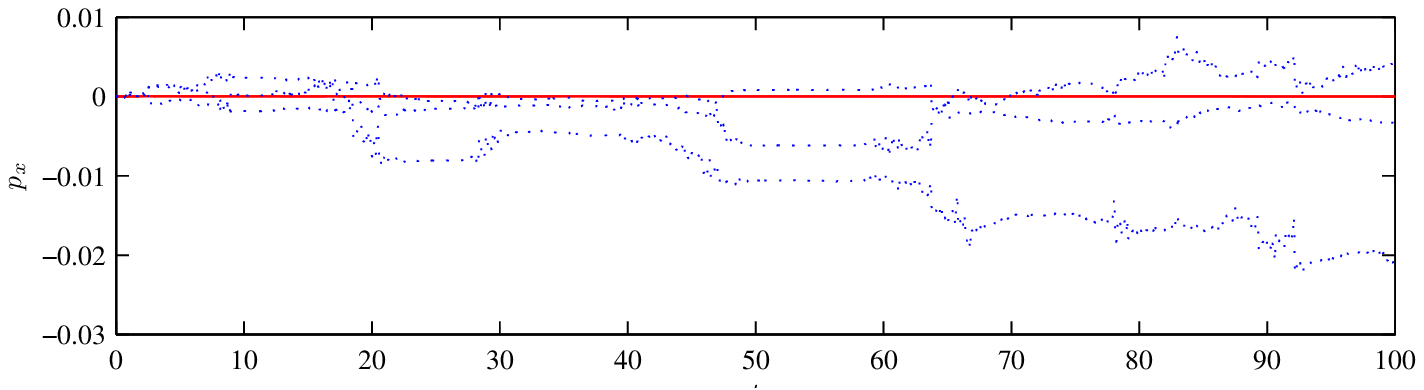}\label{fig:px}}
}
\centerline{
    \subfigure[Angular Velocity of Body 2, $\Omega_2$]{
        \includegraphics[width=0.95\columnwidth]{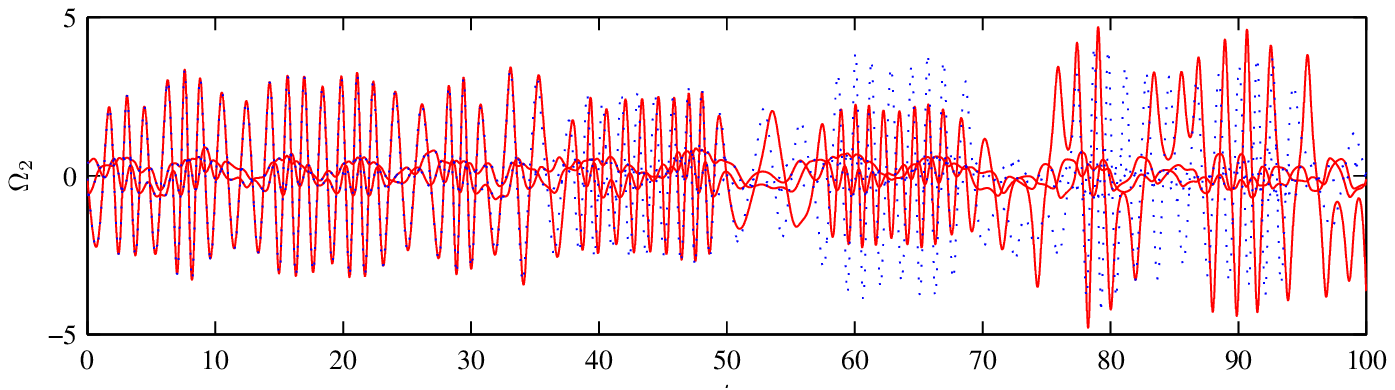}}
    \hspace*{0.05\columnwidth}
    \subfigure[Deviation of Total Angular Momentum]{
        \includegraphics[width=0.95\columnwidth]{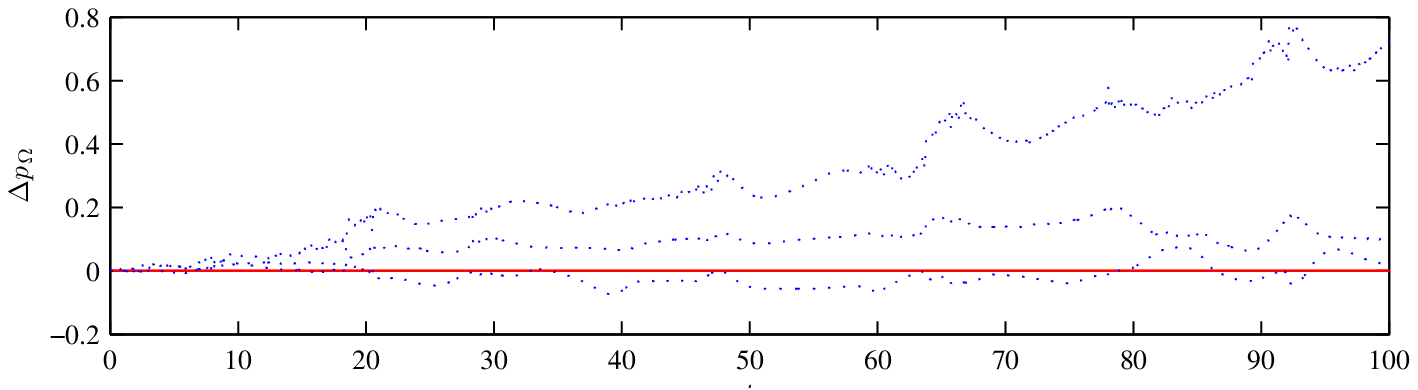}\label{fig:pw}}
}
\centerline{
    \subfigure[Velocity $\dot x$]{
        \includegraphics[width=0.95\columnwidth]{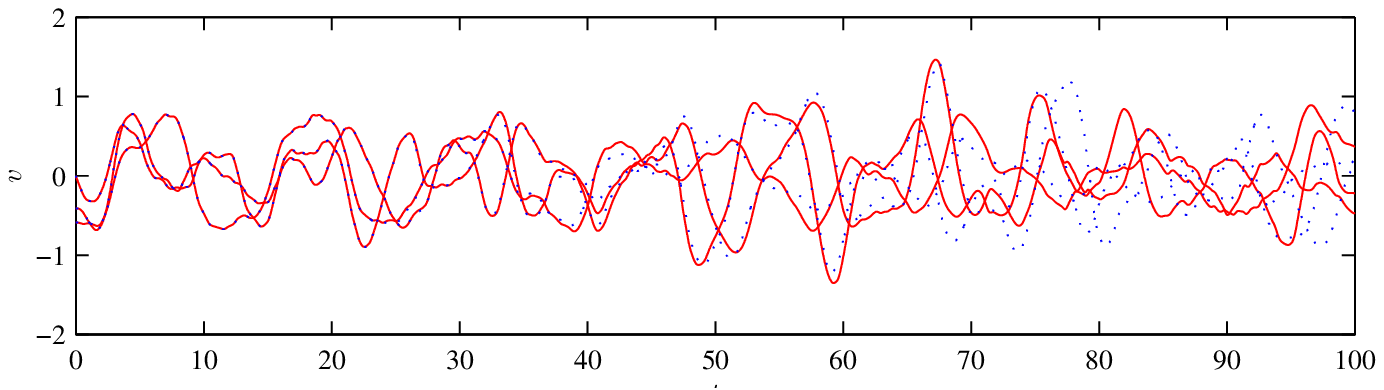}}
    \hspace*{0.05\columnwidth}        
    \subfigure[Orthogonality Error $\|I-R_i^T R_i\|,\quad i\in\{0,1,2\}$]{
        \includegraphics[width=0.95\columnwidth]{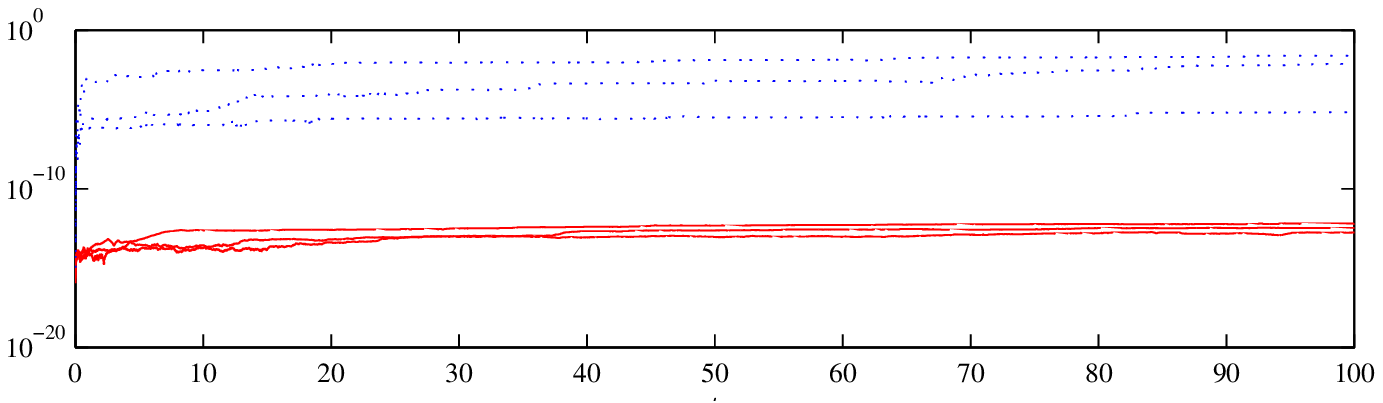}\label{fig:errR}}
}
\caption{Numerical simulation of connected rigid bodies in a perfect fluid (LGVI: red, solid, RK(4)5: blue, dotted)}\label{fig:NS}
\end{figure*}

\reffig{NS} shows the resulting angular/linear velocity responses, total energy, total linear momentum, total angular momentum deviation, and orthogonality errors of rotation matrices. The Lie group variational integrator and the Runge-Kutta method provide compatible trajectories only for a short period of time.

The computational properties of the Lie group variational integrator are as follow. As shown in \reffig{E}, the computed total energy of the Lie group variational integrator oscillates near the initial value, but there is no increasing or decreasing drift for long time periods. This is due to the fact that the numerical solutions of symplectic numerical integrators are exponentially close to the exact solution of a perturbed Hamiltonian~\cite{Hai.AoNM1994}. The value of the perturbed Hamiltonian is preserved in the discrete-time flow. The Lie group variational integrator preserves the momentum map exactly as in \reffig{px} and \ref{fig:pw}, and it also preserves the orthogonal structure of rotation matrices accurately. The orthogonality errors, measured by $\|I-R_i^T R_i\|$ for $i\in\{0,1,2\}$, are less than $10^{-13}$  in \reffig{errR}.

These show that the structure-preserving properties of the Lie group variational integrator are important for simulating the dynamics of the connected rigid bodies in a fluid accurately. A more extensive comparison study of the computational accuracy and efficiency of Lie group variational integrators can be found in~\cite{Lee.CMaDA2007}.

\section{Conclusions}

We have developed continuous-time equations of motion and a geometric numerical integrator, referred to as a Lie group variational integrator, for a system of connected rigid bodies immersed in a perfect fluid. The rigid bodies are modeled as three-dimensional ellipsoids, and each joint has three rotational degrees of freedom. This model characterizes qualitative behaviors of three-dimensional fish locomotion.

The continuous-time equations of motion provide an analytical model that is defined globally on the Lie group configuration manifold, and the Lie group variational integrator preserves the geometric features of the system, thereby yielding a reliable numerical simulation tool.

\bibliography{ACC09}

\begin{thebibliography}{10}
\providecommand{\url}[1]{#1}
\csname url@rmstyle\endcsname
\providecommand{\newblock}{\relax}
\providecommand{\bibinfo}[2]{#2}
\providecommand\BIBentrySTDinterwordspacing{\spaceskip=0pt\relax}
\providecommand\BIBentryALTinterwordstretchfactor{4}
\providecommand\BIBentryALTinterwordspacing{\spaceskip=\fontdimen2\font plus
\BIBentryALTinterwordstretchfactor\fontdimen3\font minus
  \fontdimen4\font\relax}
\providecommand\BIBforeignlanguage[2]{{%
\expandafter\ifx\csname l@#1\endcsname\relax
\typeout{** WARNING: IEEEtran.bst: No hyphenation pattern has been}%
\typeout{** loaded for the language `#1'. Using the pattern for}%
\typeout{** the default language instead.}%
\else
\language=\csname l@#1\endcsname
\fi
#2}}

\bibitem{Sfa.IJoOE1999}
M.~Sfakiotakis, D.~Lane, and J.~Davies, ``Review of fish swimming modes for
  aquatic locomotion,'' \emph{IEEE Journal of Oceanic Engineering}, vol.~24,
  no.~2, pp. 237--252, 1999.

\bibitem{Tay.PotRSoLSA1952}
G.~Taylor, ``Analysis of the swimming of long narrow animals,''
  \emph{Proceedings of the Royal Society of London. Series A}, vol. 214, no.
  1117, pp. 158--183, 1952.

\bibitem{Wu.JoFM1961}
T.~Wu, ``Swimming of a waving plate,'' \emph{Journal of Fluid Mechanics},
  vol.~10, pp. 321--344, 1961.

\bibitem{Lig.BK1975}
M.~Lighthill, \emph{Mathematical Biofluiddynamics}.\hskip 1em plus 0.5em minus
  0.4em\relax SIAM, 1975.

\bibitem{Lig.PotRSoLSB1971}
------, ``Large-amplitude elongated-body theory of fish locomotion,''
  \emph{Proceedings of the Royal Society of London. Series B}, vol. 179, pp.
  125--138, 1971.

\bibitem{Lig.ARoFM1969}
------, ``Hydromechanics of aquatic animal propulsion,'' \emph{Annual Review of
  Fluid Mechanics}, vol.~1, no.~1, pp. 413--446, 1969.

\bibitem{Nak.Pot4IOaPEC1994}
T.~Nakaoka and Y.~Toda, ``Laminar flow computation of fish-like motion wing,''
  in \emph{Proceedings of the 4th International Offshore and Polar Engineering
  Conference}, 1994, pp. 530–--538.

\bibitem{Ram.AiFM1996}
R.~Ramamurti, R.~Lohner, and W.~Snadberg, ``Computation of the unsteady-flow
  past a tuna with caudal fin oscillation,'' \emph{Advances in Fluid
  Mechanics}, vol.~9, pp. 169--178, 1996.

\bibitem{Kan.JoNS2005}
E.~Kanso, J.~Marsden, C.~Rowley, and J.~Melli-Huber, ``Locomotion of
  articulated bodies in a perfect fluid,'' \emph{Journal of Nonlinear Science},
  vol.~15, pp. 255--289, 2005.

\bibitem{Kel.1998}
S.~Kelly, ``The mechanics and control of robotic locomotion with applications
  to aquatic vehicles,'' Ph.D. dissertation, California Institute of
  Technology, 1998.

\bibitem{Rad.2003}
J.~Radford, ``Symmetry, reduction and swimming in a perfect fluid,'' Ph.D.
  dissertation, California Institute of Technology, 2003.

\bibitem{Jal.PoAViMCS1995}
J.~Jalbert, S.~Kashin, and J.~Ayers, ``A biologically-based undulatory
  lamprey-like {A}{U}{V},'' in \emph{Proceedings of Autonomous Vehicles in Mine
  Countermeasures Symposium}, 1995, pp. 39--52.

\bibitem{Bar.1996}
D.~Barrett, ``Propulsive efficiency of a flexible hull underwater vehicle,''
  Ph.D. dissertation, Massachusetts Institute of Technology, 1996.

\bibitem{Kan.PotICoDaC2005}
E.~Kanso and J.~Marsden, ``Optimal motion of an articulated body in a perfect
  fluid,'' in \emph{Proceedings of the IEEE Conference on Decision and
  Control}, 2005, pp. 2511--2516.

\bibitem{Ros.PotACC2006}
S.~Ross, ``Optimal flapping strokes for self-propulsion in a perfect fluid,''
  in \emph{Proceedings of the American Control Conference}, 2006, pp.
  4118--4122.

\bibitem{Hai.BK2006}
E.~Hairer, C.~Lubich, and G.~Wanner, \emph{Geometric Numerical Integration},
  2nd~ed., ser. Springer Series in Computational Mathematics.\hskip 1em plus
  0.5em minus 0.4em\relax Springer-Verlag, 2006, vol.~31.

\bibitem{Lee.CMaDA2007}
T.~Lee, M.~Leok, and N.~H. McClamroch, ``Lie group variational integrators for
  the full body problem in orbital mechanics,'' \emph{Celestial Mechanics and
  Dynamical Astronomy}, vol.~98, no.~2, pp. 121--144, June 2007.

\bibitem{Lee.2008}
T.~Lee, ``Computational geometric mechanics and control of rigid bodies,''
  Ph.D. dissertation, University of Michigan, 2008.

\bibitem{Hol.PD1998}
P.~Holmes, J.~Jenkins, and N.~Leonard, ``Dynamics of the {K}irchhoff equations
  {I}: Coincident centers of gravity and bouyancy,'' \emph{Physica D}, vol.
  118, pp. 311--342, 1998.

\bibitem{Lam.BK1932}
H.~Lamb, \emph{Hydrodynamics}.\hskip 1em plus 0.5em minus 0.4em\relax Cambridge
  University Press, 1932.

\bibitem{Mar.BK1999}
J.~Marsden and T.~Ratiu, \emph{Introduction to Mechanics and Symmetry},
  2nd~ed., ser. Texts in Applied Mathematics.\hskip 1em plus 0.5em minus
  0.4em\relax Springer-Verlag, 1999, vol.~17.

\bibitem{Mar.2001}
J.~Marsden and M.~West, ``Discrete mechanics and variational integrators,'' in
  \emph{Acta Numerica}.\hskip 1em plus 0.5em minus 0.4em\relax Cambridge
  University Press, 2001, vol.~10, pp. 317--514.

\bibitem{Ise.2000}
A.~Iserles, H.~Munthe-Kaas, S.~{N\o rsett}, and A.~Zanna, ``Lie-group
  methods,'' in \emph{Acta Numerica}.\hskip 1em plus 0.5em minus 0.4em\relax
  Cambridge University Press, 2000, vol.~9, pp. 215--365.

\bibitem{Hai.AoNM1994}
E.~Hairer, ``Backward analysis of numerical integrators and symplectic
  methods,'' \emph{Annals of Numerical Mathematics}, vol.~1, no. 1-4, pp.
  107--132, 1994, scientific computation and differential equations (Auckland,
  1993).

\end{thebibliography}
\bibliographystyle{IEEEtran}

\end{document}